\documentclass[12pt,amssymb]{amsart}
\usepackage[dvips]{graphicx}

\newtheorem{thm}{Theorem}[section]
\newtheorem{prop}[thm]{Proposition}
\newtheorem{cor}[thm]{Corollary}
\newtheorem{lem}[thm]{Lemma}
\newtheorem{conj}[thm]{Conjecture}
\newtheorem{exa}[thm]{Example}
\newtheorem{defn}[thm]{Definition}

\newtheorem{rem}[thm]{Remark}
\newtheorem{note}[thm]{Notation}
\newtheorem{alg}[thm]{Algorithm}

\newcommand{\ben}{\begin{enumerate}}
\newcommand{\een}{\end{enumerate}}
\newcommand{\ble}{\begin{lem}}
\newcommand{\ele}{\end{lem}}
\newcommand{\bth}{\begin{thm}}
\newcommand{\eth}{\end{thm}}
\newcommand{\bpr}{\begin{prop}}
\newcommand{\epr}{\end{prop}}
\newcommand{\bco}{\begin{cor}}
\newcommand{\eco}{\end{cor}}
\newcommand{\bcon}{\begin{conj}}
\newcommand{\econ}{\end{conj}}
\newcommand{\bde}{\begin{defn}}
\newcommand{\ede}{\end{defn}}
\newcommand{\bex}{\begin{exa}}
\newcommand{\eex}{\end{exa}}
\newcommand{\brem}{\begin{rem}}
\newcommand{\erem}{\end{rem}}
\newcommand{\bnot}{\begin{note}}
\newcommand{\enot}{\end{note}}
\newcommand{\balg}{\begin{alg}}
\newcommand{\ealg}{\end{alg}}
\usepackage{epsf}

\newcommand{\set}[1]{\left\{#1\right\}}

\newcommand{\<}{\left <}
\renewcommand{\>}{\right > }

\DeclareMathOperator{\stab}{stab}

\title[Hurwitz equivalence  is undecidable ]{The Hurwitz Equivalence Problem is Undecidable }

\author[Liberman]{E. Liberman$^1$}
\thanks{$^1$Partially supported by the Emmy Noether Research Institute for
Mathematics, Bar-Ilan University, and the Minerva Foundation, Germany, and by the Excellency Center ``Group 
Theoretic Methods in the study of Algebraic Varieties'' of the
National Science Foundation of Israel.}
\address{Department of Mathematics and Statistics, Bar-Ilan University, Ramat-Gan 52900, Israel}
\email{Liberman@macs.biu.ac.il}

\author[Teicher]{M. Teicher }
\address{Department of Mathematics and Statistics, Bar-Ilan University, Ramat-Gan 52900, Israel}
\email{teicher@macs.biu.ac.il}

\date{\today}

\begin{document}

\maketitle

\begin{abstract}
In this paper, we prove that the Hurwitz equivalence problem for
1-factorizations in $F_2 \oplus F_2$ is undecidable, and as a
consequence, the Hurwitz equivalence problem for
$\Delta^2$-factorizations in the braid groups $B_n, n\geq 5$ is
also undecidable.

\end{abstract}
\section*{Introduction}
It has long been conjectured (e.g \cite{BMT}) that the Hurwitz
equivalence problem is undecidable.  In this paper we present a
proof to this conjecture.

This problem relates to Algebraic Geometry as follows:  There is a
well defined construction (see \cite{BMT}) which attaches to any
projective curve a $\Delta^2$-factorization of the braid group,
called the \emph{braid monodromy factorization}.  This gives rise
to the definition of the braid monodromy type (BMT) of projective
curves: Two curves have the same BMT if their braid monodromy
factorizations are Hurwitz equivalent up to at most one global
conjugation.  In the same paper it is shown that two cuspidal
curves are isotopic if and only if their BMT are equal.  This
profound theorem invites the following question:  Does there exist
a finite algorithm which recognizes whether two braid monodromy
factorizations belong to the same braid monodromy type?  In order
to answer this question, we first ask a slightly simpler question:
Does there exist a finite algorithm which determines whether two
$\Delta^2$-factorizations are Hurwitz-equivalent?  In this paper
we intend to show that the answer to this question is negative.
This problem, of determining whether two factorizations are
Hurwitz-equivalent, has been discussed in many papers, but only
partial results have been achieved.

In this paper we shall find a connection between the Hurwitz-
equivalence problem and the word problem of finitely presented
groups.   The word problem is very well known, and in
(\cite{boone},\cite{novikov}) it is shown that there exist
finitely presented groups whose word problem is unsolvable.  In
(\cite{AdianI},\cite{AdianII},\cite{Rabin}) it is shown that
determining whether the word problem of a given group is solvable
is itself undecidable.  In \cite{BR} it is shown further, that
determining whether a group's word problem is solvable is
unrecognizable, and there is no uniform partial algorithm which
solves the word problem for all the finitely presented groups
whose word problem is solvable.  In \cite{Higman} Higman proves
the existence of a universal finitely presented group $K$, such
that there exists a Turing machine which, given a finitely
presented group, finds a finite subset of $K$ which generates a
subgroup isomorphic to the given group.  For more information on
the word problem and decision problems in group theory, see for
example \cite{Miller}.

This paper is organized as follows:  In chapter 1 we shall present
some well known definitions and results which we intend to
investigate.  In chapter 2 we shall study the structure of the
Hurwitz stabilizer, and the effect of the Hurwitz braid action on
direct products.  In chapter 3 we shall define a new structure and
study the connection between the Hurwitz equivalence problem and
the problem of finding a solution to the equation $H^Y = H_1$,
where $Y$ comes from a specified normal subgroup.  In chapter 4 we
present our main theorem, which connects between the word problem of
an arbitrary group and the $1$-factorizations of the braid group.
Finally, in chapter 5, we shall prove that the Hurwitz-equivalence
problem is undecidable for $F_2 \oplus F_2$, and as an important
consequence, that the Hurwitz-equivalence problem for
$\Delta^2$-factorizations in the braid group $B_n, n \geq 5$ is
undecidable.

\section{Preliminaries}

The braid group is the group
$$\left\langle  \sigma_1,\ldots,\sigma_{n-1}   \begin{array}{l}
                                \vline  \; \sigma_i\sigma_{i+1}\sigma_i=\sigma_{i+1}\sigma_i\sigma_{i+1} \\
                                  \vline\; \sigma_i\sigma_j=\sigma_j\sigma_i
                                  \hspace{1cm} if|i-j|\geq 2
                                \end{array}
                                \right \rangle$$

There is a natural homomorphism $\pi: B_n \to S_n$, defined by
mapping $\sigma_i$ to the transposition $(i,\, i+1)$.  The
permutation $\pi(b)$ of a braid $b$ is called the braid permutation.
There is a normal subgroup $PB_n$ of $B_n$, consisting of those
braids whose braid permutation is trivial.  It is shown in \cite{BU}
that the generators of $PB_n$ are $$A_{ij}  = \sigma _i^{ - 1}
\cdots\sigma _{j - 2}^{-1} (\sigma _{j - 1}^2 )\sigma _{j - 2}
\cdots\sigma _i $$, with $1 \leq i<j\leq n$.

\bth \label{artin} \cite{artin}
Let $b$ be a braid in $B_n$, and $\pi$ any braid with the same braid
permutation.  Then $b$ can be written as $\pi U_n \dots U_1$,
where each $U_i$ is generated only by generators of the form $A_{ik}$.
\eth

More information on the pure braid groups and braid combing -- the
process of rewriting a braid in this normal form -- can be found
in \cite{artin} and \cite{birman}.

In this paper we shall be particularly interested in a (right)
group action called the \emph{Hurwitz braid action}.

\bde
Let $G$ be a group.  We call a vector $F=(f_1, \dots, f_n) \in
G^n$ a $g$-factorization of length $n$ if its product
$m(F):=\prod_{i=1}^n f_i$ is equal to $g$.

Let $n$ be a natural number, and $g\in G$.  The \emph{Hurwitz braid
action} is the (right) group action of $B_n$ on the set of
$g$-factorizations of length $n$, defined as follows:  Let
$F=(f_1,\ldots,f_n)  $ be a factorization then
\begin{align*}
(f_1 ,\ldots, f_n  )\sigma _i  & = (f_1 ,\ldots, f_{i - 1} , f_{i + 1} ,f_{i + 1}^{ - 1}
f_i f_{i + 1} , f_{i + 2} ,\ldots, f_n)  \\
(f_1 ,\ldots, f_n  )\sigma _i^{ - 1} & = (f_1  ,\ldots, f_{i - 1}  , f_i f_{i + 1} f_i^{ - 1}  , f_i
, f_{i + 2}, \dots , f_n)
\end{align*}
\ede
From now on, when we refer to a group action, we mean the Hurwitz
braid action.

\bde [Hurwitz equivalence] Two factorizations are said to be
Hurwitz equivalent if they belong to the same group orbit.  We
shall denote this equivalence relation by $\cong_H$. \ede

Note than nothing crucial would change if we chose to define the
equivalence through a left group action.

\bde
Let $F=(f_1,\dots, f_n)$ be a factorization of length $n$,
and $K=(k_1,\dots, k_m)$ a factorization of length $m$.  We denote
by $F\otimes K$ the concatenation of these two factorizations,
i.e. $F\otimes K = (f_1,\dots, f_n,k_1,\dots, k_m)$.  Evidently
this is a factorization of length $m+n$.
\ede

\bde
Let $X=(g_1,\dots, g_s)$ and $Y$ be factorizations.  Then
$X^Y:=(g_1^{m(Y)}, \dots, g_s^{m(Y)})$ and
$X^{(Y^{-1})}:=(g_1^{m(Y)^{-1}}, \dots, g_s^{m(Y)^{-1}})$
\ede

We shall now give an example of how the generators $A_{ij}$ act on
a factorization.

\bex \label{ex}
Let $F$ be a factorization. Let $X$ be the first $i-1$ elements of
$F$, $a$ be the $i$-th element, $Y$ be the elements between the $i+1$-th
and the $j-1$-th elements, $b$ be the $j$-th element and $Z$ be the remaining elements of $F$, from
the $j+1$-th element and onward.  Then
\[
\begin{gathered}
  (X \otimes  a\otimes  Y \otimes b \otimes Z)A_{ij}  = (X \otimes a \otimes Y \otimes b \otimes Z)\sigma _i^{ - 1} \cdots\sigma _{j - 2}^{ - 1} (\sigma _{j - 1}^2 )\sigma _{j - 2} \cdots\sigma _i  =  \hfill \\
  (X \otimes Y^{a^{ - 1} }  \otimes a \otimes b \otimes Z)(\sigma _{j - 1}^2 )\sigma _{j - 2} \cdots\sigma _i  = (X \otimes Y^{a^{ - 1} }  \otimes b \otimes a^b  \otimes Z)(\sigma _{j - 1} )\sigma _{j - 2} \cdots\sigma _i  \hfill \\
   = (X \otimes Y^{a^{ - 1} }  \otimes a^b  \otimes b^{a^b }  \otimes Z)\sigma _{j - 2} \cdots\sigma _i  =  \hfill \\
  (X \otimes a^b  \otimes Y^{a^{ - 1} a^b  }  \otimes b^{a^{ b} }  \otimes Z) \hfill \\
\end{gathered}
\]

\[
\begin{gathered}
  (X \otimes a \otimes Y \otimes b \otimes Z)A_{ij}^{ - 1}  = (X \otimes a \otimes Y \otimes b \otimes Z)\sigma _i^{ - 1} \cdots\sigma _{j - 2}^{ - 1} (\sigma _{j - 1}^{ - 2} )\sigma _{j - 2} \cdots\sigma _i  =  \hfill \\
  (X \otimes Y^{a^{ - 1} }  \otimes a \otimes b \otimes Z)(\sigma _{j - 1}^{ - 2} )\sigma _{j - 2} \cdots\sigma _i  = (X \otimes Y^{a^{ - 1} }  \otimes b^{a^{ - 1} }  \otimes a \otimes Z)(\sigma _{j - 1} )\sigma _{j - 2} \cdots\sigma _i  \hfill \\
   = (X \otimes Y^{a^{ - 1} }  \otimes a^{(b^{a^{ - 1} } )^{ - 1} }  \otimes b^{a^{ - 1} }  \otimes Z)\sigma _{j - 2} \cdots\sigma _i  =  \hfill \\
  (X \otimes  a^{(b^{a^{ - 1} } )^{ - 1} }  \otimes Y^{ a^{ - 1}  a^ {({b^{a^{ - 1}})} ^{-1 }}}  \otimes b^{a^{ - 1} }  \otimes Z) \hfill \\
\end{gathered}
\]
\eex

\bco\label{exc}
We can see from this example that $A_{ij}$ and $A_{ij}^{-1}$
replace the $i$-th element of the factorization with a conjugate
of itself, while the other elements are either unchanged or
conjugated by an element from the normal closure of the $i$-th
element in the subgroup generated by the factorization elements.
\eco

\subsection{Recursiveness and Recursive Enumerability} \hfill

An \emph{alphabet} is a finite set $\Sigma$.  Elements of $\Sigma$ are
called \emph{symbols}.  \emph{Words} from $\Sigma$ are finite sequences of
symbols from $\Sigma$.  $\Sigma^*$ is the set of all words from
$\Sigma$.

Let $\Sigma$ be an alphabet, and $X$ a subset of $\Sigma^*$.  $X$
is called a recursively enumerable set if there exists a Turing
machine $T$ such that $T$ halts on any word from $X$, and does not
halt otherwise.  $X$ is called a recursive set if there is a
Turing machine which outputs $1$ on any input from $X$, and
outputs $0$ for inputs from $\Sigma^* \setminus X$.
Similarly, a decision problem is called solvable if the set of
elements whose answer is `yes' is recursive, and recognizable if
the set of elements whose answer is `yes' is recursively
enumerable.

\section{Properties of the Hurwitz Braid Action}

The following proposition (Proposition \ref{simple}) describes
the basic properties of the Hurwitz braid action and the Hurwitz
equivalence.  The properties are all natural, and are detailed here
for the readers convenience.

\bpr\label{simple} \hfill
   \ben
\item Let $G=(g_1,\dots ,g_n)$ and $K=(k_1,\dots, k_n)$ be
factorizations.  If $G \cong_H K$, then the group generated by the
elements of $G$ is equal to the group generated by the elements of
$K$, i.e., $\langle g_1, \dots, g_m \rangle = \langle k_1, \dots,
k_n \rangle$.

\item Let $K$ be a subgroup of $G$, and let $(g_1, \dots, g_n)$
and $(h_1, \dots, h_n)$ be two factorizations such that $g_i$ and
$h_i$ are elements from $K$.  Then $(g_1, \dots, g_n)\cong_H (h_1,
\dots, h_n)$ as factorizations in $K$ if and only if $(g_1, \dots,
g_n)\cong_H (h_1, \dots, h_n)$ as factorizations in $G$.

\item Let $(g_1, \dots, g_n)$ be a factorization and $G=\langle
g_1, \dots, g_n \rangle$.  Let $b$ be a braid with braid
permutation $\pi$.  Assume that $(g_1, \dots, g_n)b = (h_1, \dots,
h_n)$.  Then, for any $i$, $g_i$ and $h_{\pi(i)}$ are conjugate in
$G$.

\item Let $h: G_1 \to G_2$ be a group homomorphism, and let
$F=(f_1, \dots, f_n)$ be a factorization in the group $G_1$.
Denote $h(F)=(h(f_1), \dots, h(f_n))$.  Then, for any braid $b\in
B_n$, $[h(F)b] = h[(F)b]$.
\een \epr
\begin{proof}
\ben
\item \label{one} We first prove the proposition for $\sigma_i$.
Let $F=(f_1, \dots, f_n)$ be a factorization, and $G_1=\<
f_1, \dots, f_m \>$ be the group generated by its elements.
Similarly, Let $G_2$ be the group generated by the elements of
$F\sigma_i = (f_1\dots f_{i-1}, f_{i+1}, f_{i+1}^{-1}f_i f_{i+1},
f_{i+2}, \dots, f_n)$, i.e. $G_2=\langle f_1\dots f_{i-1},
f_{i+1}, f_{i+1}^{-1}f_i f_{i+1}, f_{i+2}, \dots, f_n \rangle$.
All the elements besides the $i$-th and $i+1$-th are equal, so in
order to prove equality, it suffices to show that $f_i,
f_{i+1} \in G_2$ and that $f_{i+1}, f_{i+1}^{-1}f_i f_{i+1} \in
G_1$.  Evidently $f_{i+1}$ belongs to both groups.  Now, since
$f_i$ and $f_{i+1}$ are in $G_1$, so is $f_{i+1}^{-1}f_i f_{i+1}$.
Similarly, since $f_{i+1}$ and $f_{i+1}^{-1}f_i f_{i+1}$ belong to
$G_2$, so does $f_i = (f_{i+1})(f_{i+1}^{-1}f_i
f_{i+1})(f_{i+1}^{-1})$.  This completes the proof on the
generators of the braid group, and the proof can easily be
completed by induction.

\item From (\ref{one}) we see that if $F=(f_1, \dots,
f_n)$ is a factorization such that $f_i \in K$, then all the
elements in the factorizations belonging to the orbit of $F$ also belong
to $K$. Now by hypothesis $K$ is a subgroup of $G$,
so the action of $\sigma_i$ on $F$ is the same, whether we
consider $F$ a factorization in $G$ or in $K$.  Therefore,
regardless of our point of view, the group orbit of $F$ remains
the same.  Since by definition $F\cong_H F_1$ if $F_1$ is in the
group orbit as $F$, we have completed the proof.

\item We first give a proof for the generators $\sigma_i$.  All
the elements in the factorization other than the $i$-th and $(i+1)$-th elements
are left unchanged by the braid action, and the braid permutation
$(i,\, i+1)$ also leaves them fixed, so the theorem holds
trivially.  It remains to check the $i$-th and $(i+1)$-th
elements.  $\sigma_i$ maps $(f_i, f_{i+1})$ to $(f_{i+1},
f_{i+1}^{-1} f_i f_{i+1})$, and we see clearly that the $(i+1)$-th
element of the first factorization is equal (and hence conjugate)
to the $i$-th element in the second one, and that the $i$-th element of
the first factorization is conjugate to the $(i+1)$-th element
of the second one.  The proof can now be completed by a straightforward
induction, bearing in mind than from $(\ref{one})$
we know that every factorization that is Hurwitz equivalent to our
factorization has its elements in $G$ too.

\item Again, we shall start by proving the proposition for
$\sigma_i$.  We know that $\sigma_i$ affects only the $i$-th and
$(i+1)$-th elements of the factorization, so we need only consider
them.  We know that the $i$-th and $(i+1)$-th elements of $h(F)$
are $h(f_i)$ and $h(f_{i+1})$ respectively, so $(h(f_i),
h(f_{i+1}))\sigma_i = (h(f_{i+1}),
h(f_{i+1}^{-1})h(f_i)h(f_{i+1}))$.  since $h$ is a group
homomorphism, this is equal to $(h(f_{i+1}), h(f_{i+1}^{-1} f_i
f_{i+1}))$, which are precisely the $i$-th and $(i+1)$-th elements of
$h[(F)\sigma_i]$.  Once again, the proof is easily completed by induction.
 \een
\end{proof}

\ble\label{stab of function}
Let $h: G_1 \to G_2$ be a homomorphism, and $F$ a factorization of
$G_1$.  Then $\stab(F) \subseteq \stab(h(F))$.
\ele
\begin{proof}
Let $b\in \stab(F)$, which means that $Fb=F$.  From Lemma \ref{simple} we have that
$[h(F)]b = h[(F)b]$.  Putting both together we obtain
$[h(F)]b=h(F)$, which implies that $b\in \stab(h(F))$.
\end{proof}

\ble\label{lemma stab}
Let $V=(x_1, \dots, x_r)$ be a factorization such that $x_i \neq
1$ for all $i$, and let $S$ be the stabilizer of $V$.  Let $\varphi:
B_r \to B_{n+r}$ be defined by $\sigma_i \mapsto \sigma_{i+n}$.
Let
$$V_1  = (\underbrace {1 ,\ldots, 1}_n , x_1 ,\ldots, x_r)$$
Then, the stabilizer of $V_1$ is $S_1=\< \set{\sigma_i \mid
i<n }, \set{A_{jk} \mid j<n}, \varphi(S) \>$.

Furthermore, we can write any element of the stabilizer in the
form $CBA$, where $A\in \< \sigma_i \mid i<n \>$, $B\in \< A_{jk} \mid
j<n \>$ and $C\in \varphi(S)$.
\ele
\begin{proof}
$S_1 \subseteq \stab(V_1)$. \hfill\\
We shall show this by considering the three types of braids in $S_1$.
\begin{itemize}
\item $\sigma_i, i<n$.  $\sigma_i$ affects only the $i$-th and
$(i+1)$-th elements, which in our case are $(1,1)$.  The operation
of $\sigma_i$ on these gives $(1, 1^{-1}\cdot 1 \cdot 1)=(1,1)$.

\item $A_{jk}, j<n$.  From example \ref{ex} we know that if $V$ is a
factorization with the same notation as in the example, then
$(V)A_{jk}  = (X \otimes a^{ba^{ - 1} }  \otimes Y^{a^{ - 1} a^{ba^{ - 1} } }  \otimes b^{a^{ - 1} }  \otimes Z)$.
However, in our case $a=1$, so we get
$$(X \otimes 1^{b*1^{ - 1} }  \otimes Y^{1^{ - 1} 1^{b*1^{ - 1} } }
\otimes b^{1^{ - 1} }  \otimes Z) = (X \otimes 1 \otimes Y \otimes b
\otimes Z) = V
$$

\item $\varphi(S)$.  The elements of $\varphi(S)$ act only on the last
$r$ elements of the factorizations, and do so the same way $S$
acts on $V$.  Therefore, $\varphi(S)\subseteq \stab(V_1)$.
\end{itemize}

$\stab(V_1) \subseteq S_1$.\hfill \\
Let $b$ is a braid in the
stabilizer.  From Theorem \ref{artin} we know that we can write $b$ as
$b= \pi F_{r+n-1}\cdots F_1$, where $\pi$ is an arbitrary braid
with the same braid permutation as $b$, and the $F_i$ is generated by
braids of the form $A_{ik}$, where $i$ is fixed and $i < k$. By hypothesis
we know that $x_i\neq 1$ for all $i$, so the $x_i$\,s cannot be conjugate to $1$.
Therefore, $b$'s braid permutation must permute the elements $1,
\dots, n$ among themselves and the elements $n+1, \dots, n+r$
among themselves.  This means that $\pi$ can be factored as $\pi_1
\pi_2$, where $\pi_1$ is generated by braids from $\set{
\sigma_i \mid i>n}$ and $\pi_2$ is generated by braids from $\set{
\sigma_i \mid i<n}$.

Now, since the pure braid group is a normal subgroup and\\
$F_{r+n-1}\cdots F_1 \in PB_n$, then $\pi_2 F_{r+n-1}\cdots F_1 = P
\pi_2$, for some \\$P\in PB_n$.  This means that $P$ can be written
as $F'_{r+n-1}\cdots F'_1$, so we have $\pi_2 F_{r+n-1}\cdots F_1 =
F'_{r+n-1}\cdots F'_1 \pi_2$.  Therefore, $b=\pi _1 \pi _2 F_{r + n
- 1} \cdots F_1  = \pi _1 (F'_{r + n - 1} \cdots F_1 ')\pi _2 $.

We know that $b$ is in the stabilizer, so $(V_1 )b = V_1$, which
means $(V_1)\pi _1 (F'_{r + n - 1} \cdots F_1 ')\pi _2  = V_1$, or
in other words\\, $(V_1 )(\pi _1 F'_{r + n - 1} \cdots F'_n )(F'_{n
- 1} \cdots F'_1 \pi _2 ) = V_1 $. The braid permutation of $(\pi _1
F'_{r + n - 1} \cdots F'_n )$ acts only on the first $n$ elements,
hence the first $n$ elements are conjugate to
 $1$, which means they are equal to $1$. The last $r$ element are conjugated among
  themselves
 and therefore they are not equal to $1$.  Therefore, we can implement the part of the lemma which we proved
 above on $(V_1 )(\pi _1 F'_{r + n - 1} \cdots F'_n)$, and conclude that $F'_{n - 1} \cdots F'_1 \pi _2$ is in the stabilizer of  $(V_1 )(\pi _1 F'_{r + n - 1} \cdots F'_n )$.

This in turn implies that $(V_1 )(\pi _1 F'_{r + n - 1} \cdots F'_n ) = V_1$,
so therefore $(\pi _1 F'_{r + n - 1} \cdots F'_n )$ is in the stabilizer of $V_1$.
Note also that it is generated by braids of the form $\{ \sigma _i \mid i > n\} $,
 and therefore it is in $\varphi (S)$.  In summary, we have shown that the braid $b$ can be written as
 $[\pi _1 F'_{r + n - 1} \cdots F'_n ][F'_{n - 1} \cdots F'_1]\pi _2$,
 where $ \pi _2\in \<\{ \sigma _i |i < n\} \>,[F'_{n - 1} \cdots F'_1] \in \<\{ A_{jk} |j < n\} \>,[\pi _1 F'_{r + n - 1} \cdots F'_n ] \in
\varphi (S)$, which completes the proof.
\end{proof}

\bde
Let $G=A \oplus B = \set{ (a, b) \mid a\in A, b\in B }$.  Let
$V=(v_1, \dots, v_n)$ be a factorization in $A$ and $W=(w_1,
\dots, w_2)$ a factorization in $B$. We denote by $V\oplus W$ the
factorization of $A\oplus B$ whose $i$-th element is $(v_i,
w_i)$.  In other words, $V\oplus W=((v_1,w_1),\ldots,(v_n,w_n))$.
\ede

\begin{lem}\label{s}
Let $b$ be a braid, then $(V\oplus W)b = (V)b \oplus (W)b$.
\end{lem}
\begin{proof}
First we shall prove the lemma for $\sigma_i$.  We know that
$\sigma_i$ affects only the $i$-th and $(i+1)$-th elements of the
factorization.  Let us write the $i$-th element of $V\oplus W$ as
$(a_1, b_1)$ and the $(i+1)$-th element as $(a_2, b_2)$.  Now,\\
$((a_1,b_1),(a_2,b_2))\sigma_1 = ((a_2, b_2),
(a_2,b_2)^{-1}(a_1,b_1)(a_2,b_2))$\\
$ = ((a_2,b_2), (a_2^{-1}a_1a_2, b_2^{-1}b_1 b_2))$, which are
precisely the $i$-th and $(i+1)$-th elements of $(V)\sigma_i \oplus
(W)\sigma_i$.  A similar argument holds for $\sigma_i^{-1}$.  The
proof for arbitrary braids follows by induction.
\end{proof}

\begin{lem}\label{lem+}
$V \oplus W \cong_H V_1\oplus W$ if and only if there exists a
braid $b\in \stab(W)$ such that $(V)b = V_1$.
\end{lem}
\begin{proof}
By definition, $V\oplus W \cong_H V_1 \oplus W$ if and only if
there is a braid $b$ such that $(V\oplus W) b = V_1 \oplus W$.
However, from lemma \ref{s} we know that $(V\oplus W)b = (V)b \oplus
(W)b$, so out requirement becomes $(V)b \oplus (W)b = V_1 \oplus W$.
This implies that $(W)b=W$, i.e. that $b\in \stab(W)$, and that
$(V)b = V_1$.
\end{proof}

\section{On the Set of Factorizations}\
In this section we construct a new structure on the set of factorizations
and study its properties.

\bde \label{p} Let $G = G_1 \oplus G_2$. Let $R_1,
\ldots,R_n, W_1,\ldots ,W_m, H \in G_1$ and $X_1,\ldots,X_{m+2} \in G_2$. Let $R$
be the factorization $(R_1,\ldots,R_n)$, $W$ the factorization $(W_1,\ldots,W_m)$ and $X$
the factorization $(X_1,\ldots,X_{m+2})$.  We define  $P_X (R,W,H)$ to be the following
 factorization on $G$:
$$P_X (R,W,H)=((R_1 ,1),\ldots, (R_n ,1) , (W_1 ,X_1 ),\ldots, (W_m ,X_m ) , (H^{-1},X_{m + 1}
), (H ,X_{m + 2} ))$$
\ede

We will now present two theorems which give necessary and
sufficient conditions for the equivalence of the Hurwitz equivalence
problem and the problem of finding an element $Y$ inside a normal
subgroup such that $H^Y=H_1$.

\bth \label{mth1} .
Let $R$ be the factorization $(R_1,\ldots,R_n)$, $W$ be the
factorization $(W_1,\ldots,W_m)$ and $X$ be the factorization
$(X_1,\ldots,X_{m+2})$.  Let $H$ and $H_1$ be elements of $G_1$.

Suppose that $R_1 \cdots R_n  = 1$.  Let $N$ be the normal closure of
$\{ R_1 ,\ldots,R_n\}$ inside $\< R_1 ,\ldots,R_n ,W_1 ,\ldots,W_m
\>$, and assume that $N$ contains a word $Y$ such that $H^Y  = H_1$.
Then $P_X (R,W,H) \cong _H P_X (R,W,H_1)$. \eth

 In order to prove Theorem \ref{mth1}, we first prove two lemmas.

\ble\label{l1}
  Let $R$ be the factorization $(R_1,\ldots,R_n)$,  $W$ be the factorization $(W_1,\ldots,W_m)$
  and $X$ be the factorization $(X_1,\ldots,X_{m+2})$.  Let $H$ be an element from $G_1$.
  Suppose also that $R_1\cdots R_n  = 1$.
  Then for any j:  $P_X (R,W,H) \cong _H P_X (R^{W_j },W ,H)$
\ele

\begin{proof}
We know from Lemma \ref{lem+}  that  $P_X (R,W,H) \cong _H P_X
(R^{W_j } ,W,H)$ if and only if there exists a braid $b \in B_{m+2}$
in the stabilizer of  $ (\underbrace {1 , \ldots, 1}_n , X_1,
\ldots, X_{m + 2}) $,
 such  that $(R_1,\ldots, R_n  , W_1  ,\ldots, W_m  , H^{ - 1}  , H)b$
is equal to \\ ($R_1^{W_j},\ldots, R_n^{W_j}   , W_1  , W_m , H^{-1}  , H $).

With this in mind, define the  braid\\
 $b_j  = (\sigma _{j - 1} \cdots \sigma _n )(\sigma _{n -
1} \cdots\sigma _1 )(\sigma _1 \cdots \sigma _{n - 1} )(\sigma _{j - 1} \cdots \sigma _n
)^{ - 1} $ for  $j>n$.

We begin by showing that  $b_j$ is in the stabilizer of
$ (\underbrace {1, \ldots, 1}_n , X_1,\ldots,   X_{m + 2}) $.\\
Let $I=(\underbrace{1,\ldots,1}_n)$, $Y=(X_1,\ldots,X_j)$  and $Z=(X_{j+1},\ldots,X_{m+2})$
and so:
$$\begin{array}{l}
(\underbrace {1 , \ldots, 1}_n , X_1,\ldots,   X_{m + 2})b_j=
 (I \otimes Y \otimes X_j  \otimes Z)b_j  \\
  = (I \otimes Y \otimes X_j  \otimes Z)(\sigma _{j - 1} \cdots\sigma _n )(\sigma _{n - 1} \cdots\sigma _1 )(\sigma _1 \cdots\sigma _{n - 1} )(\sigma _{j - 1} \cdots\sigma _n )^{ - 1}  \\
  = (I \otimes X_j  \otimes Y^{X_j }  \otimes Z)(\sigma _{n - 1} \cdots\sigma _1 )(\sigma _1 \cdots\sigma _{n - 1} )(\sigma _{j - 1} \cdots\sigma _n )^{ - 1}  \\
  = (X_j  \otimes I^{X_j }  \otimes Y^{X_j }  \otimes Z)(\sigma _1 \cdots\sigma _{n - 1} )(\sigma _{j - 1} \cdots\sigma _n )^{ - 1}  \\
 \end{array}$$
Since $I$ contains only $1$'s, every element of  $I^{X_j } $
 is $1^{m(X_j )}  = 1$ and therefore $I^{X_j }  = I$.  Therefore, we get
\begin{align*}
(I \otimes Y \otimes X_j  \otimes Z)b_j & = (X_j  \otimes I \otimes Y^{X_j }  \otimes Z)(\sigma _1 \cdots\sigma _{n - 1} )(\sigma _{j - 1} \cdots\sigma _n )^{ - 1}  \\
 & = (I \otimes (X_j )^I  \otimes Y^{X_j }  \otimes Z)(\sigma _{j - 1} \cdots\sigma _n )^{ - 1}
 \end{align*}
Moreover, since $I$ is a factorization of $1$'s,  $m(I)^{-1}  = 1$ and therefore we get
\begin{align*}
 & = (I \otimes X_j  \otimes Y^{X_j }  \otimes Z)(\sigma _{j - 1} \cdots\sigma _n )^{ - 1}  \\
 &= (I \otimes Y \otimes X_j  \otimes Z) \\
\end{align*}
It remains to be shown that $(R_1,\ldots, R_n  , W_1  , W_m  , H^{ - 1} , H)b_j$ is equal to\\
$(R_1^{W_j}  , R_n^{W_j} , W_1  , W_m  , H^{ - 1}  , H) $. Using the more concise notation
$R=(R_1,\ldots,R_n)$, $X=(W_1,\ldots,W_{j-1})$, $Y=(W_{j+1},\ldots,W{m},H,H^{-1})$, and starting off from the
left hand side of the equation, we have:
\begin{align*}
(R_1,\ldots, &R_n, W_1  , W_m  , H^{ - 1} , H)b_j =\\
  & = (R \otimes X \otimes W_j  \otimes Y)b_j  \\
  & = (R \otimes X \otimes W_j  \otimes Y)(\sigma _{j - 1} \cdots\sigma _n )(\sigma _{n - 1} \cdots\sigma _1 )(\sigma _1 \cdots\sigma _{n - 1} )(\sigma _{j - 1} \cdots\sigma _n )^{ - 1}  \\
  & = (R \otimes W_j  \otimes (X)^{W_j }  \otimes Y)(\sigma _{n - 1} \cdots\sigma _1 )(\sigma _1 \cdots\sigma _{n - 1} )(\sigma _{j - 1} \cdots\sigma _n )^{ - 1}  \\
  & = (W_j  \otimes R^{W_j }  \otimes (X)^{W_j }  \otimes Y)(\sigma _1 \cdots\sigma _{n - 1} )(\sigma _{j - 1} \cdots\sigma _n )^{ - 1}  \\
  & = (R^{W_j }  \otimes (W_j )^{R^{W_j } }  \otimes (X)^{W_j }  \otimes Y)(\sigma _{j - 1} \cdots\sigma _n )^{ - 1}  \\
\end{align*}
Now, we have by hypothesis that $m(R)=1$, and therefore also $m(R^{W_j } ) = 1$, so we can continue:
\begin{align*}
 & = (R^{W_j }  \otimes W_j  \otimes (X)^{W_j }  \otimes Y)(\sigma _{j - 1} \cdots\sigma _n )^{ - 1}  \\
 & = (R^{W_j }  \otimes ((X)^{W_j } )^{W_j^{ - 1} }  \otimes W_j  \otimes Y) \\
 & = (R^{W_j }  \otimes (X)^{W_j W_j^{ - 1} }  \otimes W_j  \otimes Y) \\
 & = (R^{W_j }  \otimes X \otimes W_j  \otimes Y) \\
 \end{align*}
\end{proof}

\ble \label{l2} Let $R$ be
the factorization $(R_1,\ldots,R_n)$, $W$ be the factorization $(W_1,\ldots,W_m)$ and $X$ be the
factorization $(X_1,\ldots,X_{m+2})$.  Let $H$ be an element from $G_1$.
Then for any $j$: $P_X (R,W,H) \cong _H P_X (R,W,H^{R_j } )$

\ele

\begin{proof}
The proof of this Lemma is very similar to that of the previous
lemma.  We know from lemma \ref{lem+} that $P_X (R,W,H) \cong _H P_X
(R,W,H^{R_j } )$
 if and only if there exist a braid $b \in B_{n+m+2}$ in the stabilizer of
 $(\underbrace {1 , \ldots, 1}_n , X_1,\ldots,   X_{m + 2}) $
 such that
 $$(R_1 ,\ldots, R_n  , W_1  ,\ldots, W_m  , H^{ - 1}  , H)b =(R_1  ,\ldots, R_n  , W_1  ,\ldots, W_m  , (H^{ - 1} )^{Rj} , (H)^{R_j} )$$

Let us define the  braid:$$ c_j  = (\sigma _j^{ - 1} \cdots\sigma _{m + n}^{ - 1}
)(\sigma _{m + n + 1} \sigma _{m + n + 2}^2 \sigma _{m + n + 1} )(\sigma _{m + n}
\cdots\sigma _j ) , j \le n$$
and show that it fulfills these conditions. We begin by showing that  $c_j $ is in the stabilizer of
 $(\underbrace {1 , \ldots, 1}_n , X_1,\ldots,   X_{m + 2}) $.  Using the shorter notation
 $I=(\underbrace{1,\ldots,1}_j)$, $Z=(\underbrace{1, \ldots, 1}_{n-(j+1)},X_1,\ldots,X_m)$ and $T=(X_{m+1},X_{m+2})$,
we see that $c_j$ is indeed in the stabilizer, because:
\begin{align*}
(\underbrace {1 , \ldots, 1}_n , &X_1,\ldots,   X_{m + 2})c_j \\
  & = (I \otimes 1 \otimes Z \otimes T)c_j  \\
  & = (I \otimes 1 \otimes Z \otimes T)(\sigma _j^{ - 1} \cdots\sigma _{m + n}^{ - 1} )(\sigma _{m + n + 1} \sigma _{m + n + 2}^2 \sigma _{m + n + 1} )(\sigma _{m + n} \cdots\sigma _j ) \\
  & = (I \otimes (Z)^{1^{ - 1} }  \otimes 1 \otimes T)(\sigma _{m + n + 1} \sigma _{m + n + 2}^2 \sigma _{m + n + 1} )(\sigma _{m + n} \cdots\sigma _j ) \\
  & = (I \otimes Z \otimes 1 \otimes T)(\sigma _{m + n + 1} \sigma _{m + n + 2}^2 \sigma _{m + n + 1} )(\sigma _{m + n} \cdots\sigma _j ) \\
  & = (I \otimes Z \otimes T \otimes 1^T )(\sigma _{m + n + 2} \sigma _{m + n + 1} )(\sigma _{m + n} \cdots\sigma _j ) \\
  & = (I \otimes Z \otimes T \otimes 1)(\sigma _{m + n + 2} \sigma _{m + n + 1} )(\sigma _{m + n} \cdots\sigma _j ) \\
  & = (I \otimes Z \otimes 1 \otimes T^1 )(\sigma _{m + n} \cdots\sigma _j ) \\
  & = (I \otimes Z \otimes 1 \otimes T)(\sigma _{m + n} \cdots\sigma _j ) \\
  & = (I \otimes 1 \otimes Z^1  \otimes T) \\
  & = (I \otimes 1 \otimes Z \otimes T)
 \end{align*}

To conclude the proof, it remains to be proved that
$$(R_1 ,\ldots, R_n  , W_1  ,\ldots, W_m  , H^{ - 1} , H)b =
(R_1  ,\ldots, R_n  , W_1  ,\ldots, W_m  , (H^{ - 1} )^{Rj} , (H)^{R_j } )$$
Again, using the shorthand notation $X=(R_1,\ldots R_{j-1})$, $Y=(R_{j+1},\ldots ,R_{n},W_1,\ldots,W_m)$,
$H'=(H^{-1}, H)$, we have:
\begin{align*} (
 R_1 ,\ldots, & R_n  , W_1  ,\ldots, W_m  , H^{ - 1} , H)c_j = \\
 & = (X \otimes R_j  \otimes Y \otimes H')c_j  \\
 & = (X \otimes R_j  \otimes Y \otimes H')(\sigma _j^{ - 1} \cdots\sigma _{m + n}^{ - 1} )(\sigma _{m + n + 1} \sigma _{m + n + 2}^2 \sigma _{m + n + 1} )(\sigma _{m + n} \cdots\sigma _j ) \\
 & = (X \otimes Y^{(R_j )^{ - 1} }  \otimes R_j  \otimes H')(\sigma _{m + n + 1} \sigma _{m + n + 2}^2 \sigma _{m + n + 1} )(\sigma _{m + n} \cdots\sigma _j ) \\
 & = (X \otimes Y^{(R_j )^{ - 1} }  \otimes H' \otimes R_j ^{H'} )(\sigma _{m + n + 2} \sigma _{m + n + 1} )(\sigma _{m + n} \cdots\sigma _j )
\end{align*}
We continue by noting that $m(H')=m((H^{-1}, H )) = 1$, so we get:
\begin{align*}
 & (X \otimes Y^{(R_j )^{ - 1} }  \otimes H' \otimes R_j )(\sigma _{m + n + 2} \sigma _{m + n + 1} )(\sigma _{m + n} \cdots\sigma _j ) \\
  & = (X \otimes Y^{(R_j )^{ - 1} }  \otimes R_j  \otimes (H')^{R_j } )(\sigma _{m + n} \cdots\sigma _j ) \\
  & = (X \otimes R_j  \otimes \left( {Y^{(R_j )^{ - 1} } } \right)^{R_j }  \otimes (H')^{R_j } ) \\
  & = (X \otimes R_j  \otimes Y^{(R_j )^{ - 1} R_j }  \otimes (H')^{R_j } ) \\
  & = (X \otimes R_j  \otimes Y \otimes (H')^{R_j } ) \\
 \end{align*}
 \end{proof}

Now that the two Lemmas have been proven,  we are ready to tackle Theorem \ref{mth1}.
\begin{proof}{(Theorem \ref{mth1})} \hfill\\
We begin  by proving that if  $x \in \<W_1 ,\ldots,W_n \>$
 then
$P_X (R,W,H) \cong _H P_X (R^x ,W,H)$ for any $H$.
Since $x \in \<W_1, \ldots,W_n \>$,
 it can be written as $x = a_1 \cdots a_\ell $
where $a_j  \in \{ W_1 ,\ldots,W_m ,W_1^{ - 1} ,\ldots,W_m^{ - 1} \} $.

We proceed by induction on $\ell$. If  $\ell=0$ then the theorem holds trivially. Next, suppose that
the theorem holds for $\ell$, and we shall prove that it holds for $\ell+1$. Since $x$ is the product of $\ell+1$
elements from $\{ W_1 ,\ldots,W_m ,W_1^{ - 1} ,\ldots,W_m^{ - 1} \} $, we can write $x$ as
$x_\ell h_{\ell + 1} $, where $x_\ell $ is the product of $\ell$ elements,
  and $h_{\ell + 1}  \in \{ W_1 ,\ldots,W_m ,W_1^{ - 1} ,\ldots,W_m^{ - 1} \} $.
We know from the  induction hypothesis that $P_X (R,W,H) \cong _H P_X (R^{x_\ell },W ,H)$.
 Since $R_1  \cdots R_n  = 1$, also  the product  $(R_1 )^{x_\ell }  \cdots (R_n )^{x_\ell }  = 1$,
 and therefore,  if  $h_{\ell + 1}  = W_j $,
 we conclude from Lemma \ref{l1} that $P_X (R^{x_\ell },W, H) \cong _H P_X(((R^{x_\ell } )^{W_j },W  ,H) = P_X (R^{x_\ell W_j } ,W,H) = P_X (R^x ,W,H)$

If on the other hand $h_{\ell + 1}  = W_j^{ - 1} $,
 then we know that  $(R_1 )^x  \cdots (R_n )^x  = 1$,
 and from Lemma \ref{l1} we get that
 $P_X (R^x ,W,H) \cong _H P_X ((R^x )^{W_j },W ,H) = P_X (R^{xW_j } ,W,H) = P_X (R^{x_\ell W_j^{ - 1} W_j } ,W,H) = P_X (R^{x_\ell } ,W,H)$

Now that we have this in hand, we can complete the proof of the
theorem. Let $Y \in N $.   we are going to prove that $P_X (R,W,H)
\cong _H P_X (R,W,H^Y
 )$.
Since  $Y \in  N $ we can write $Y=(x_1 )^{y_1 } (x_2 )^{y_2 }
,\ldots,(x_l )^{y_l } $
 where $x_i  \in \{ R_1 ,\ldots,R_n ,R_1^{ - 1} ,\ldots,R_n^{ - 1} \} $
 and $y_i  \in \<W_1 ,\ldots,W_n\>$.

We prove the theorem by induction on $l$. If $l=0$,  the  theorem holds trivially.
Assume that the theorem holds for $l$ and we shall prove it for $l+1$. Indeed, suppose that
$Y = (x_1 )^{y_1 } (x_2 )^{y_2 } \cdots(x_{l + 1} )^{y_{l + 1} }$.  $Y$ can be factored
thus: $Y=Y_l({x_{l+1}})^{y_{l+1}}$, where $Y_l  = (x_1 )^{y_1 } (x_2 )^{y_2 } \cdots(x_l )^{y_l } $.
By the induction hypothesis $P(R,W,H) \cong _H P(R,W,H^{Y_l } )$.

Suppose now that that  $x_{l + 1}  = (R_j )^{y_{l + 1} }$.
We have already proved that $P_X (R,W,H^{Y_l } ) \cong _H P_X(R^{y_{l+1} } ,W,H^{Y_l} )$,
and from lemma \ref{l2} we know that
$
  P_X(R^{y_{l+1} } ,W,H^{Y_l} ) \cong_H P_X (R^{y_{l + 1} } ,W,(H^{Y_l } )^{R_j^{y_{l + 1} } } )$\\
$ = P_X (R^{y_{l + 1} } ,W,H^{Y_l R_j^{y_{l + 1} } } ) = P_X
(R^{y_{l + 1} } ,W,H^Y )$.
 Since $y_{l+1}  \in \<W_1 ,\ldots,W_n \>$,  by what we have already proved we get that  $P_X
(R^{y_{l + 1} } ,W,H^Y ) \cong _H P_X (R,W,H^Y )$

On the other hand, if  $x_{l + 1}  = (R_j^{ - 1} )^{y_{l + 1} } $,
 since $y_{l+1}\in \<W_1,\ldots,W_m\>$, by what we initially proved we get that  $P_X (R,W,H^Y ) \cong _H P_X (R^{y_{l+1} } ,W,H^Y )$
and from lemma \ref{l2} we get that\\
$ P_X (R^{y_{l + 1} } ,W,H^Y ) \cong _H P_X (R^{y_{l + 1} } ,W,(H^Y
)^{(R_j )^{y_{l + 1} } } ) = P_X (R^{y_{l + 1} } ,W,H^{YR_j^{y_{l +
1}} }   ) $ \\ $= P_X (R^{y_{l + 1} } ,W,H^{Y_l (R_j^{ - 1} )^{y_{l
+ 1} }R_j^{y_{l + 1} }
 }  ) = P_X (R^{y_{l + 1} } ,W,H^{Y_l } )
$, and since $y_{l+1}  \in \<W_1 ,\ldots,W_n \>$, $P_X (R^{y_{l + 1} },W ,H^{Y_l } \cong_H
P_X (R,W,H^{Y_l } )$
 \end{proof}

  \bth\label{mth2}
 Let $ G:=G_1 \oplus G_2 $, and
\begin{align*}
    F_1 & := ((x_1,1),\ldots,(x_n,1),(y_1,z_1),\ldots,(y_m,z_m)) \\
    F_2 & :=((x_1',1),\ldots,(x_n',1),(y_1',z_1),\ldots,(y_m',z_m)) \\
    K & :=\<x1,\ldots,x_n,y_1,\ldots,y_m\>
\end{align*}
Let $N$ be the normal closure of $\{x_1,\ldots,x_n\}$ inside $K$.

\ben

\item If the following conditions hold:
\ben
\item $F_1\cong_H F_2$,
\item $z_i \ne 1$, for all $i$,
\een
then there exists a permutation $\pi\in S_n$ such that for any $i$, $1\leq i\leq n$, $x_i'$
is conjugated to $x_{\pi(i)}$ in $K$ (in particular $x_i' \in N$).
\item If in addition the braid-action stabilizer of $(z_1,\ldots,z_m)$ is contained in the stabilizer
of $(y_1,\ldots,y_m)$
 then for any $i$,
 $1\leq i\leq m$,  there exists an element $m_i \in N$ such that $y_i'={y_i}^{m_i}$
 (Note that both the stabilizers are subgroups of $B_m$ and
 therefore the requirement that one be contained in the other is well defined.)
 \een
\eth
\begin {proof}
\ben
\item
 Since $F_1\cong_H F_2$, there is a permutation $\pi_1 \in S_{m+n}$ such that the $r$-th element in $F_2$
is conjugated to the $\pi_1(r)$-th element in $F_1$.
  Let $i$ be a number, $1\leq i\leq n$.  We know that $(x_i',1)$ is conjugated in $G$ to the $\pi_1(i)$-th element
  of the factorization $F_1$, say $(a,b)$. This means that $x_i'$ is conjugated to $a$ in
$G_1$ and $1$ is conjugated to $b$ in $G_2$, which implies that $b=1$.  Now, since $z_j\neq 1$ for
all $j$, $(x_i',1)$ must be conjugated to $(x_l,1)$ for some
$l$, $1\leq l\leq n$.  This in turn implies that $x_i'$ is conjugated to $x_l$ in $K$, and that $\pi_1$ induces a
permutation on the first $n$ elements.

\item
We know from Lemma \ref{lem+} that $F_1 \cong_H F_2$ if and only if there exists a braid $b$
in the stabilizer of $(\underbrace {1 ,\ldots, 1}_n,z_1,\ldots,z_m)$ such that
$(x_1,\ldots,x_n,y_1,\ldots,y_m)b$ is equal to
 $(x_1',\ldots,x_n',y_1',\ldots,y_m')$.  From Lemma \ref{lemma stab} we know that if $S$ is the stabilizer of
   $(z_1,\ldots,z_m)$,  then every element in  the stabilizer of $(\underbrace {1 ,\ldots, 1}_n,z_1,\ldots,z_m)$  can
be written  as $CBA$ where
$A \in \<\{ \sigma _i |i < n\} \>$, $B \in \<\{ A_{jk} |j < n,j\leq k\leq n+m\} \>$, $C \in \varphi
(S)$. In particular, we can write $b=CBA$, with $A, B, C$ as above.  Since $C \in \varphi (S)$,
 and since we know that the stabilizer of $(z_1,\ldots,z_m)$
 is contained in the stabilizer of $(y_1,\ldots,y_m)$, then $C$ must be in the stabilizer of $(x_1,\ldots,x_n,y_1,\ldots,y_m)$,
so $(x_1,\ldots,x_n,y_1,\ldots,y_m)BA = (x_1',\ldots,x_n',y_1',\ldots,,y_m')$.  In order to prove the theorem,  we must prove
that for any braid of the form $BA$, with $A \in \<\{ \sigma _i \mid i < n\} \>$ and $B \in \<\{
A_{jk} \mid j < n,j\leq k\leq n+m\} \>$, there exists a factorization $X$ of length $n$ and
$p_1,\ldots,p_m \in N$ such that $(F_1)BA=X\otimes (y_1^{p_1},\ldots,y_m^{p_m})$.

We first prove this fact for $B$, by induction on the length of $B$ when written as a product
of the generators $ \{ A_{jk} \mid j < n,j\leq k\leq n+m\}$.  If the length is $0$, then
the theorem is trivial. Assume now that the theorem is true for $\ell$ and we will
prove it for $\ell+1$.  Let $B=B_\ell b_{\ell+1}$, where $B_\ell$ is a braid of length
$\ell$ and $b_{\ell+1}$ is a braid of the form $A_{ir}^\varepsilon$, where $i<n$, $i\leq
r\leq m+n$ and $\varepsilon \in \{ \pm 1\}$. From the induction hypotheses, we conclude that there
exists
a factorization of length $n$ and  $p_1,\ldots,p_m \in N $,  such that
$(x_1,\ldots,x_n,y_1,\ldots,y_m)B_\ell=X\otimes (y_1^{p_1},\ldots,y_k^{p_k})$. From Corollary \ref{exc}
we know that $A_{ir}^\varepsilon$ acts on every element except
the $i$-th element of $X\otimes (y_1^{p_1},\ldots,y_k^{p_k})$
 by  conjugating it by an element of the normal closure of the $i$-th element inside $K$.
 In particular, since $i<n$, it acts on all the elements beyond the $n$-th one by
 conjugating them by elements from normal closure of the $i$-th element.
 From the first part of the theorem we know that the $i$-th element is in $N$,  which implies
 that its normal closure is also contained in $N$, so there is a factorization $X'$ of length $n$, and $P_1,\ldots,P_m\in N$
  such that $X\otimes (y_1^{p_1},\ldots,y_k^{p_k})A_{ir}^\varepsilon=X'\otimes({(y_1^{p_1})}^{P_1},\ldots,{(y_k^{p_k})}^{P_k})
  =X'\otimes(y_1^{p_1P_1},\ldots,y_k^{p_kP_k})$. We conclude the induction by noting that for any $j$, $p_j$ and $P_j$ are
  in $N$, whence also $p_jP_j\in N$.

To complete the proof, we must prove the fact for $A$.
However, since $A$ is generated by braids of the form $\{ \sigma _i  \mid i <
n\}$, all the elements beyond the $n$-th element remain fixed.
\een

\end {proof}

\bco\label{mco}
 Let $R$ be the factorization $(R_1,\ldots,R_n)$, $W$ be the factorization $(W_1,\ldots,W_m)$
  and $X$ be the factorization $(X_1,\ldots,X_{m+2})$.  Let $H$ and $H_1$ be elements from $G_1$.
Assume that the following conditions hold: \ben
\item $R_1\cdots R_n=1$.
\item $H \in \<R_1 ,\ldots,R_n ,W_1 ,\ldots,W_m \>$.
\item For every $i$,  $X_i  \ne 1$.

\item The braid-action stabilizer of $(X_1  ,\ldots, X_{m + 2} )$
 is contained in the stabilizer of $(W_1  ,\ldots , W_m  , H^{-1} , H )$. (Note that both of the stabilizers are subgroups of $B_{m+2}$,
 therefore the question of whether one is contained in the other is well defined.)
 \een
 Then  the normal subgroup of $\<R_1 ,\ldots,R_n ,W_1 ,\ldots,W_m \>$
 generated by $\{R_1 ,\ldots,R_n\} $
 contains a word Y such that  $H^Y  = H_1 $ if and only if $P_X (R,W,H) \cong _H P_X (R,W,H_1 )$.
\eco
\begin{proof}
 The \emph{if} part follows trivially from Theorem \ref{mth1}.
For the \emph{only if} direction, we use Theorem \ref{mth2}, which implies that if
$P_X (R,W,H) \cong _H P_X (R,W,H_1 )$, then there is an element $Y$ in the normal closure
 of $\<R_1,\ldots,R_n,W_1,\ldots,W_m,H^{-1},H\>$, generated
 by $\{R_1,\ldots,R_n\}$ such that $H^Y=H_1$. Since $H \in \<R_1 ,\ldots,R_n ,W_1 \ldots,W_m
 \>$, $\<R_1,\ldots,R_n,W_1,\ldots,W_m,H^{-1},H\>=\<R_1,\ldots,R_n,W_1,\ldots,W_m\>$ and $Y$ is indeed in
 the normal closure of $\{R_1,\ldots,R_n\}$ in $\<R_1,\ldots,R_n,W_1,\ldots,W_m\>$

\end{proof}
\section{the functions $FTL_n$}
We shall now show that the word problem and the Hurwitz equivalence
problem are equal in some sense. We first present some
well known notation on the free group. A word is said to be
\emph{reduced} if there are no two adjacent letters of the form
$x_ix_i^{-1}$ or $x_i^{-1}x_i$.  A word is \emph{cyclic reduced}
if it is reduced and the first and the last letters are not
inverses. We say that a word $A$ has a \emph{root} if there is a
word $w$ and $n \in \mathbb{Z}, n\neq \pm 1$, (or equivalently, $n\geq
2$) such that $A=w^n$.

We begin by defining the sets:
$$S_1:=\left\{(F_n,V_m)  \begin{array}{l}
                       \vline \mbox{$F_n$ is a free group generated by $z_1,\ldots,z_n$ for some $n$, $1\leq n <\infty$}   \\
                        \vline V_m=(v_1,\ldots,v_m)\in (F_n^*)^m, 0\leq m<\infty
                     \end{array}
                     \right \}$$

$$NS_1:=\{(F_n,V_m,a) \mid (F_n,V_m)\in S_1, a \in {F_n}^*\}$$

Next, we define a function $FTL_{s1}$ from $NS_1$ to the set of
$1$-factorization of the group $F_2\oplus F_2$.

\begin{defn}[{\emph{$FTL_{s1}$}}] \hfill\\
Let $(F_n,V_m,a)\in NS_1$. Let $F_X$ be the free group with
$n+1$ generators, generated by $X_1,\ldots,X_{n+1}$, where these symbols
are taken to be distinct of the symbols in $(F_n,V_m,a)$.

Note that for every $r_1,r_2 \in \mathbb{N}$ $F_{r_1}\oplus
F_{r_2} \subseteq F_2 \oplus F_2$, so we can chose a family of
injections $i_{(r_1,r_2)}:F_{r_1}\oplus F_{r_2}\rightarrow F_2
\oplus F_2$. We define $FTL_{s1}(F_n,V_m,a):=i_{(n,n+1)}P_X(
Vf,W,a)$ where $Vf:=(v_1,\ldots, v_m,(v_1\cdots v_m)^{-1})$,
$W:=(z_1,\ldots,z_n,(z_1\cdots z_n)^{-1})$ and
$X:=(X_1,\ldots,X_n,(X_1\cdots X_n)^{-1},X^{-1}_{n+1},X_{n+1})$.
We can easily see that $FTL_{s1}(F_n,V_m,a)$ is a $1$-factorization.
\end{defn}

\bth\label{s1} Let $(F_n,V_m) \in S_1$.  Let $a,b$ be
words from $F_n$ and $N$ be the normal closure of $V_m$ in $F_n$.  Then
$FTL_{s1}(F_n,V_m,a)\cong_HFTL_{s1}(F_n,V_m,b)$ if and only if
there exists an element $Y \in N$  such that $a^Y=_{F_n}b$ \eth

\begin{proof}
It is sufficient to prove that $P_X(Vf,W,a)\cong_H P_X(Vf,W,b)$ if
and only if there exists $Y$ in the normal closure of $V_m$ in $F_n$
such that $a^Y=_{F_n}b$. This will follow from Corollary \ref{mco},
once we've confirmed that all the conditions of the corollary hold.
But indeed,

\ben

\item The product of the elements of $Vf$ is $1$

\item Every element in $X$ is different from $1$.

\item The word $a$ is in $F_n$.

\item  The stabilizer of X is contained in the stabilizer of $W
\otimes a^{-1} \otimes a$, for any $a$. This is true because we
can construct an homomorphism $h$ from $F_X$ to $G_W$ by
$X_i\mapsto z_i$ for $1\leq i\leq n$ and $X_{n+1}\mapsto a$. Since
$h$ is a homomorphism, by Lemma \ref{stab of function} the
stabilizer of $X$  is contained in stabilizer of $W \otimes a^{-1}
\otimes a=h(X)$

\een

We can therefore conclude from Corollary \ref{mco} that
$P_X(Vf,W,a) \cong _H P_X(Vf,W,b)$ if and only if there exists an
element $Y$ such that $a^Y=_{F_n}b$, where
$Y$ is in the normal closure of $\{v_1,\ldots,v_m,(v_1 \cdots v_m)^{-1}\}$ in
$\<v_1,\ldots,v_m,(v_1 \cdots v_m)^{-1},z_1,\ldots,z_n,(z_1\cdots
z_n)^{-1}\>$, which is equal to the normal closure of
$\{v_1,\ldots,v_m\}$ in $\<z_1,\ldots,z_n\>$ (i.e. $N$)
\end{proof}
\begin{lem} \label{Z} \hfill
\begin{enumerate}
\item\label{Z1} If $a$ and $b$ are elements of a free group $F$, such that $a^m$ and
$b^n$ commute, then $a$ and $b$ must be powers of a common element $c$.
\item If $a$ is an element with no root in the free group $F$, then
its centralizer is $Z(a)=\{a^k \mid k\in\mathbb{Z}\}$
\end{enumerate}
\end{lem}
\begin{proof}
\begin{enumerate}
\item See \cite{Group}.
\item Obviously $a$ and $a^k$ commute for any $k\in \mathbb{Z}$ so
what remains to be shown is that if $b\in Z(a)$, then
$b=a^k$ for some $k \in \mathbb{Z}$. But indeed, if $b\in Z(a)$, then by
definition $a$ and $b$ commute, which by (\ref{Z1}) means that they
are both powers of a common element $c$, i.e. $a=c^m$
and $b=c^n$.  However, since $a$ has no root, we must conclude that $m=\pm1$, which implies
$c=a^{\mp 1} $, so in summary we have $b=c^m={( a^{\mp 1})}^m=a^{\mp m}$, which concludes the proof.
\end{enumerate}
\end{proof}
 \ble
\label{root} Let $F_n, (n \ge 2)$ be a free group, and let $N$ be a
non trivial normal subgroup. There exists an element $A$ in $N$ such
that $A$ has no root \ele
\begin{proof}
 Let $y \ne e$ be
 some element in $N$. Since $N$ is normal we may assume that $y$ is cyclic-reduced.
 We shall separate our proof into two cases: \\
\emph{Case 1:}  $y = (x_i )^n $. \\
 Without loss of generality we can take $y = (x_1 )^n$.
 Now define $A = x_2^{ - 1} (x_1 )^n x_2 (x_1 )^n$.  Note that $A$ is cyclic
 reduced.  We claim that $A$ has no root.  Indeed, assume by negation that there is a cyclic reduce word
 $w$ such that for some $m\geq 2$, $w^m  = A$.
 $w$ must contain  $x_2^{-1}$
 which means that $w^m $
 must contain $x_2^{-1} $ $m$ times,
 which is not the case.\\
\emph{Case 2:} $y$ contains at least two different generators. \\
In this case, we can
write our word (perhaps after cyclic permutation) in such a way
that it begins and ends with different generators. Without loss of
generality we can write the word as $y = x_1 Z x_2^\varepsilon $
 where $\varepsilon  \in \{ 1, - 1\} $.
 Denote the length of  $y$ by $l$, and define
  $A = (x_2^\varepsilon x_1   )^{-(l + 1)} x_1 Zx_2^\varepsilon  (x_2^\varepsilon  x_1  )^{  (l + 1)} x_1
  Zx_2^\varepsilon$.  Again, $A$ is cyclic reduced and we claim that it has no root.
  Suppose, to the contrary, that there is a word $w$ such that for some $m\geq2$
  $w^m  = A$.  Since our word start with $(x_2^\varepsilon x_1 )^{-(l + 1)}$ and also $x_1 $ appears in $A$,
 also $w$ must start with $(x_2^\varepsilon x_1)^{-(l + 1)}$.
 In addition, since our word contain $(x_2^\varepsilon  x_1  )^{  (l + 1)} $
 and since $w$ start with $x_1^{-1} $,  $w$ must contain  $(x_2^\varepsilon  x_1  )^{  (l + 1)} $.
Summing up, we see that the length of $w$ is at least $4l+4$, so the length of $w^m$ is at least $m(4l+4) > 6l+4$.
However, the length of  $A$ is only $6l+4$, and we have reached a contradiction.
\end{proof}

Let us define the sets:
$$S_2:=\{(F_n,V_m) \mid (F_n,V_m,a) \in NS,n\geq
2,\exists  v \in V_m \mbox{ such that } v\neq_{F_n} 1 \}. $$
$$NS_2:=\{(F_n,V_m,a) \mid (F_n,V_m,a) \in NS_1,(F_n,V_m)\in S_2 \}. $$
 We shall now construct a function from $NS_2$ to the set of
$1$-factorizations in the group $F_2 \oplus F_2$, which we denote
$FTL_{s2}$.

\bde [\emph{$FTL_{s2}$}]  \hfill\\
Let $(F_n,V_m,a)\in NS_2$, and let $N$ be
the normal closure of $V_m$ in $F_n$.  $V_m$ contains a word
which represents a non-trivial element in $F_n$, so $N$ cannot be
trivial. Since $n\geq 2$ and N is non-trivial,
    from Lemma \ref{root} we conclude that there is a word $H$ in $N$ which has no root.
    We fix such an $H$ for $F_n$ and $V_m$,
    and define $FTL_{s2}(F_n,V_m,a)=FTL_{s1}(F_n,V_m,H^a)$

\ede

\bth\label{s2} Let $(F_n,V_m)\in S_2$.  Let
 $a,b$ be words in $F_n$ and $N$ be the normal closure of $V_m$ in $F_n$.
 Then $FTL_{s2}(F_n,V_m,a)\cong FTL_{s2}(F_n,V_m,b)$ if and only if $a^{-1}b \in
 N$.
 \eth
 \begin{proof}
By definition  $FTL_{s2}(F_n,V_m,a)\cong_H FTL_{s2}(F_n,V_m,b)$ if and only if
 $FTL_{s1}(F_n,V_m,H^a)\cong_H FTL_{s1}(F_n,V_m,H^ b)$.
 From Theorem \ref{s1} we know that  $FTL_{s1}(F_n,V_m,H^a)\cong_H FTL_{s1}(F_n,V_m,H^b)$
 if and only if $\exists Y \in N$ such that $(H^a)^Y=H^b$. Now, $Y_1=a^{-1}b$ is readily seen to be a solution
 of the equation $(H^a)^Y=H^b$, so the general solution is $XY_1$, where $X$ commutes with $H^a$.
However, since $H$ has no root, neither does $H^a$, so from Lemma \ref{Z}
we know that only $(H^a)^k$ commute with $H^a$, and therefore $Y$ must
be equal to $(H^a)^k a^{-1}b$ for some $k\in \mathbb{Z}$. Therefore,
$FTL_{s2}(F_n,V_m,a)\cong FTL_{s2}(F_n,V_m,b)$
 if and only if $(H^a)^ka^{-1}b \in N$ for some $k\in \mathbb{Z}$.  However, $H \in N$ and $N$ is normal in $F_n$, so $(H^a)^k\in N$
 for any $k \in \mathbb{Z}$.  From this we conclude that $FTL_{s2}(F_n,V_m,a)\cong FTL_{s2}(F_n,V_m,b)$ if and only if
 $a^{-1}b \in N$
 \end{proof}

Let us now define the set $GS$ of all groups with at least
two generators and a non trivial relation, and the set
$$GSA:=\{(G,a) \mid G\in GS,a\in G^*\}$$
We shall now construct a function from $GSA$ to
 the set of  $1$-factorizations of $F_2\oplus F_2$, which we denote $FTL_B$.

 \bde{\emph{$FTL_{B}$}\\} Let  $G:= \< W_1,\ldots,W_n \mid R_1,\ldots,R_m \>
 $, and let $a$ be a word in $G$. Let $F_n$ be the free group with the
  generators $z_1,\ldots,z_n$, and let $\varphi$ be the rewriting
  function which maps words in $G$ to words in $F_n$ by replacing $W_i$ with $z_i$ and $W_i^{-1}$
  with $z_i^{-1}$. We define
  $FTL_B(G,a)=FTL_{s2}(F_n,(\varphi(R_1),...,\varphi(R_m)),\varphi(a))$.
This function is well defined, because $G$ has at least two
generators and a non-trivial relation.

\ede

 \bth\label{finitely presented}
 Let $G \in GA$  and let $a,b$ be words in $G$.
 Then $a=_Gb$ if and only if $FTL_B(G,a)\cong_H FTL_B(G,b)$.
\eth

\begin{proof}
 Let  $N$ be the normal closure of $\{\varphi(R_1),...,\varphi(R_m)\}$ in $F_n$.
We know that $a=_G b$ if and only if $\varphi(a^{-1}b)\in N$.
By Theorem \ref{s2} we know that this happens if and only if
$$FTL_{s2}(F_n,(\varphi(R_1),...,\varphi(R_m)),\varphi(a))\cong_HFTL_{s2}(F_n,(\varphi(R_1),...,\varphi(R_m)),\varphi(a))$$
which, by definition, is equivalent to $FTL_B(G,a)\cong_HFTL_B(G,b)$

\end{proof}

 \bth \label{em}  $F_2  \oplus F_2  \subseteq B_5 \subseteq  B_n, (n>5) $ \eth
 \begin{proof}
 \cite{Makanina}
 \end{proof}

\bde[\emph{$FTL_n$}]\label{ftln}\hfill\\
From Theorem \ref{em} we know that there exists a family of injections
$ib_n:F_2\oplus F_2\rightarrow B_n$. We define the functions
$FTL_n$ from $GS$ to the set of $1$-factorizations in $B_n$  by $FTL_n(G,a)=ib_n(FTL_B(G,a))$\ede

 \bco\label{bnc}  Let $G \in GS$, and let $a,b$ be words in $G$.
Then $a=_Gb$ if and only if  $FTL_n(G,a)\cong_HFTL_n(G,b)$\eco

Note that groups with a single generator or with no relations can
be changed into the desired form  by adding a generator and a relation
which equates it to the identity.  We can also notice that these groups
are the cyclic groups and the free groups, which have uniform
algorithms for solving their word problems. Therefore, we may assume that our
function is defined on all finitely presented groups.

\section{The Main Results}
The following theorem is our main undecidability result.

\begin{thm}\hfill
\begin{enumerate}
\item  The Hurwitz equivalence problem for $1$-factorizations and
$\Delta^2$-factorizations in $B_n\,\, (n\geq 5)$ is undecidable.

\item For any $n\geq 5$ there exist a $1$-factorization of $B_n$,
which we denote by $u$, such that the problem of determining
whether a $1$-factorization $v$ is Hurwitz-equivalent to $u$ is
undecidable.

\item For any $n\geq 5$ there exist a $\Delta^2$-factorization of
$B_n$, denoted by $u$, such that the problem of determining
whether a $\Delta^2$-factorization $v$ is Hurwitz-equivalent to
$u$ is undecidable.
\end{enumerate}
\end{thm}

\begin{proof} \hfill
\begin{enumerate}
\item Follows from (2) and (3)

\item Let $G$ be some finitely presented group with an undecidable
word problem.  Let $FTL_n$ be the function into
$\{$1-factorizations in $B_n\}$ which we constructed in
\ref{ftln}. Define $u:=FTL_n(G,1)$, a $1$-factorization of $B_n,
(n\geq 5)$. Suppose now that there exists an algorithm which
decides whether a given $1$-factorization $v$ is  Hurwitz
equivalent to $u$ or not. In particular, this algorithm would
decide whether $FTL_n(G,1) \cong_H FTL_n(G,a)$ for any word $a$ in
$G$. By \ref{bnc}, this is the same as deciding whether $a =_G 1$,
which contradicts our choice of $G$ as having an undecidable word
problem.

\item Let $u$ and $v$ be two $1$-factorizations. Since $\Delta^2$
is central, $u \cong_H v$ if and only if $\Delta^2 \otimes u
\cong_H \Delta^2 \otimes v$.  However, $\Delta^2 \otimes u$ and
$\Delta^2 \otimes v$ are $\Delta^2 $-factorizations, which proves
the theorem.
\end{enumerate}
\end{proof}

\begin{defn}
Let $u$ be a factorization and $X$ a set of factorizations.  $u$
is said to be \emph{compatible} with $X$ if there exists a finite
algorithm which, when given any $v \in X$, decides whether $u
\cong_H v$ or not.
\end{defn}

The following theorem gives unrecognizability results.

\begin{thm} \hfill
\begin{enumerate}
\item There exists a recursive subset $X$ of the
$\Delta^2$-factorizations of $B_n \,\, n\geq 5$, such that if we take
the subset
$$X_1=\set{u \in X \mid u \text{ is compatible with }X},$$
then:
\begin{enumerate}
\item the Hurwitz-equivalence problem on $X_1$ is undecidable.

\item $X_1$ is not recursively enumerable.

\item The problem of determining whether a
$\Delta^2$-factorization is compatible with $X$ is unrecognizable.
\end{enumerate}

\item The problem of determining whether a recursive subset of the
$\Delta^2$-factorizations of $B_n$, $n\geq5$ has a decidable
Hurwitz-equivalence problem is unrecognizable.

\item There exists a recursive subset $X$ of the
$\Delta^2$-factorizations such that

\begin{enumerate}
\item the problem of determining whether a recursive subset of $X$
has solvable Hurwitz equivalence problem is unrecognizable.

\item For any $u$ in $X$, the problem of determining whether a
recursive subset of $X$ is compatible with $u$ is unrecognizable.

\end{enumerate}
\end{enumerate}
\end{thm}
\begin{proof}
\begin{enumerate}

\item We know that the set of presentations of finitely presented
groups is a recursive set.  This means that there exists an
injective total recursive function $h$ from the presentations of
finitely presented groups into $\mathbb{N}$.  We now define the
set
$$X:=\left\{\Delta^{2(h(G)+1)}\otimes \Delta ^{-2h(G)}\otimes
FTL_n(G,A) \vert \begin{array}{l}\vline \mbox{$G$ is a finitely
presented group} \\ \vline \mbox{$A$ is a word in $G$}
\end{array}\right\}$$

Let $u \in X$, so $u$ can be written as $\Delta^{2(h(G)+1)}\otimes
\Delta ^{-2h(G)}\otimes FTL_n(G,A)$ where $G$ is a finitely
presented group and $A$ is a word in $G$.  Since for any finitely
presented group $G'$ and any word $A$ from $G'$ there is no
element of the form $\Delta^{2k}, k\in \mathbb{Z}$, and since $h$
is injective and $\Delta^{2k}$ is in the center for any $k \in
\mathbb{Z}$, the only elements that can be Hurwitz-equivalent to
$u$ are $\set{\Delta^{2(h(G)+1)}\otimes \Delta ^{-2h(G)}\otimes
FTL_n(G,A) \mid A \text{ is a word in } G}$.  This set has
decidable Hurwitz equivalence problem if and only if $G$'s word
problem is solvable, so the subset containing all the
factorizations for which there exists an algorithm which
determines whether they are Hurwitz-equivalent to some
factorization from $X$ is $\{\Delta ^{2(h(G)+1)}\otimes \Delta
^{-2h(G)}\otimes FTL_n(G,A)\mid G$ is a finitely presented group
with solvable word problem and $A$ is a word in $G\}$. We can now
prove our theorem.
\begin{enumerate}

\item Suppose, by contradiction, that there exists an algorithm
which solves the Hurwitz-equivalence in this set.  This would
imply that there is an uniform algorithm to solve all word
problems for groups with solvable word problems, which we know is
impossible.

\item\label{2} Assume by negation that $X_1$ is recursively
enumerable.  This would mean that there is an algorithm which
halts on an element in $X_1$ if and only the element is compatible
with $X$. Now, letting $G$ be a finitely generated group, we can
recognize whether $G$'s word problem is solvable by checking
whether $FTL_n(G,1)$ is compatible with $X$.  This is again
impossible.

\item Suppose that the set of compatible factorizations is
recursively enumerated. We know that $X$ is a recursive set, so
their intersection, which is $X_1$, is recursively enumerable,
contradicting \ref{2}.

\end{enumerate}

\item Let $G$ be a finitely presented group and $W$ the set of all
its words.  Determining whether the set $\{\Delta^2 \otimes
FTL_n(G,A) \mid A \in W\}$ has decidable Hurwitz-equivalence problem is
equivalent to determining the problem of whether $G$'s word
problem is unsolvable, which is known to be an unrecognizable
problem.

\item We know \cite{Higman} that there exists a universal finitely
presented group $K$ and a Turing machine $T$ such that for any
finitely presented group $G$, $T(G)$ is a finite set of words in
$K$ such that $\langle T(G) \rangle \simeq G$. We denote
$X:=\{FTL_n(K,A) \mid \mbox{$A$ is a word in $K$}\}$.

\begin{enumerate}
\item Let $G$ be a finitely presented group. To determine whether
the recursive set \\ $\{\Delta^2 \otimes FTL_n(K,A) \mid \mbox{$A$
is a word generate by $T(G)$ and their inverses}\}$ has solvable
Hurwitz problem is equivalent to determining whether $G$ has
solvable word problem, which is an unrecognizable problem.

\item Let $u \in X$, so we can write $u$ as $\Delta^2 \otimes
FTL_n(K,B)$ where $B$ is a word in $K$. Let $G$ be a finitely
presented group.  In order to recognize whether $G$ has solvable
word problem we can ask whether for any word $C \in \langle
T(G)\rangle$ we can tell whether $C=_K 1$, which is equivalent to
asking whether $CB=_K B$. So in conclusion, if we consider the
recursive set which is also subset of $X$:\\
$\{\Delta^2 \otimes FTL_n(K,AB) \mid \mbox{$A$ is a word generated
by $T(G)$ and their inverses}\}$,  recognizing whether this set is
compatible with $u=\Delta^2 \otimes FTL_n(K,B)$ is equivalent to
recognizing whether $G$ has solvable word problem, which is an
unrecognizable problem.

\end{enumerate}
\end{enumerate}
\end{proof}
\brem All the theorems in this chapter will hold if we take instead
$B_n$ any group including $F_2\oplus F_2$ and instead $\Delta^2$ an
element from the center. We give now a proof for one theorem but the
rest of the proofs are similar  \erem \bth let $G$ be a group
including $F_2\oplus F_2$ and $g\in Z(G)$ then the    Hurwitz
equivalence problem for $g-factorization$ is undecidable \eth
\begin{proof}
since $F_2\oplus F_2 \subseteq G$ there is a injective homomorphism
$i:F_2\oplus F_2 \rightarrow G$ we shell define a function as in the
braid group from $GSA$ to the $g-factorization$ by $g\otimes
i(FTL_{s2}(K,a)$ where K is a finitely  presented  group and a word
from this group. We could see that the hurwitz equivalence problem
for $g-factorization$ is undecidable  while the proof is very
similar from now on to the case of the braid group.
\end{proof}

\end{document}